\documentclass[12pt, a4paper]{article}
\usepackage[latin1]{inputenc}
\usepackage{amsmath}
\usepackage[english]{babel}
\usepackage{indentfirst}
\usepackage{amstext}
\usepackage{amsfonts}
\usepackage{textcomp}
\usepackage{amssymb}
\usepackage[dvips]{graphicx}
\usepackage{setspace}
\usepackage[all]{xy}

\newcommand {\demo}{\hskip -0.6cm{\bf Proof.  }}
\newcommand {\fim}{\hfill{$\square$}\vskip 1pc}

\newcommand {\R}{\mathbb{R}}

\newcommand {\C}{\mathbb{C}}

\newcommand{\supp}{\text{supp}}

\newcommand{\id}{\mathrm{id}}

\newtheorem{teorema}{Theorem}[section]

\newtheorem{obs}[teorema]{Remark}

\newtheorem{definicao}[teorema]{Definition}
\newtheorem{proposicao}[teorema]{Proposition}
\newtheorem{exemplo}[teorema]{Example}

\begin{document}

\onehalfspace

\title{Faithful representations of graph algebras via branching systems}

\date{11 Dec 2014}

\author{Daniel Gon\c{c}alves\footnote{Partially supported by CNPq}, Hui Li\footnote{Corresponding author. Supported by a Lift-off Fellowship from Australian Mathematical Society and by Research Center for Operator Algebras of East China Normal University} \ and Danilo Royer}

\maketitle

AMS 2000 MSC: 46L05, 37A55

Keywords: C*-algebra, graph algebra, Leavitt path algebra, branching system, representation.

\begin{abstract}
We continue to investigate branching systems of directed graphs and their connections with graph algebras. We give a sufficient condition under which the representation induced from a branching system of a directed graph is faithful and construct a large class of branching systems that satisfy this condition. We finish the paper by providing a proof of the converse of the Cuntz-Krieger uniqueness theorem for graph algebras by means of branching systems.
\end{abstract}

\section{Introduction}

Directed graphs are combinatorial objects that appear in numerous situations throughout all mathematical subjects.  In particular, graph C*-algebras were introduced about two decades ago, see \cite{flr, KPRR}, as generalizations of Cuntz-Krieger algebras and more recently, see \cite{AAP, AAP1}, algebraic analogues of graph C*-algebras, called Leavitt path algebras, were introduced. Both graph C*-algebras and Leavitt path algebras (which here forth we just call graph algebras) have been the focus of intense research in the last few years, one of the main reasons for this is the fact that many combinatorial properties of a directed graph characterize properties of the associated algebra and vice versa.

It is natural to consider the relations between Leavitt path algebras and graph C*-algebras. Actually, the study of these relations was one of the main goals of the meeting "Bridges between graph C*-algebras and Leavitt path algebras" which was held in April 2013 at BIRS, Canada. Among the motivating aspects for the study of these relations is the fact that many results of graph C*-algebras have Leavitt path algebras versions and vice versa. For example, the graph-theoretic conditions under which the C*-algebra $C^*(E)$ of a directed graph $E$ is simple (finite-dimensional, AF, simple purely infinite, respectively) are precisely the same as the graph-theoretic conditions under which the Leavitt path algebra $L_K(E)$ is simple (finite-dimensional, ultramatricial, simple purely infinite, respectively). However, there is no prescription on how to obtain a result in one setting from a similar result in the other setting. In fact their proofs are often completely independent, which leads to the development of new methods.

In the previous work \cite{GR4, MR2903145, GR3, GR, GR2}, motivated by the connection between wavelet theory and representations of the Cuntz-Krieger algebra, see \cite{Bratteli}, the study of representations of graph algebras via branching systems has been initiated and developed. Branching systems arise in many areas of mathematics such as the Perron-Frobenius operator from ergodic theory (see \cite{MR2903145, GR}). In \cite{GR2} it was shown that for a large class of directed graphs every representation of a graph C*-algebra is unitarily equivalent to a representation induced from a branching system (a similar result for Leavitt path algebras was shown in \cite{GR3}). Furthermore, in the Leavitt path algebra context and in case of row-finite directed graphs without sinks, in \cite{GR3} a sufficient condition over a branching system, to guarantee faithfulness of the induced representation, was given. In this paper, we find an analogous condition over branching systems of an arbitrary graph and prove, by completely different means, that the representation of the graph C*-algebra from a branching system satisfying such condition is faithful. 

Finally, we take advantage of branching systems techniques to give an alternative proof of the converse of the Cuntz-Krieger uniqueness theorem for  graph algebras. In the context of graph C*-algebras this result can be derived from a more general result by Katsura (see \cite[Theorem~6.14]{Katsura:ETDS06}). The advantage of our proof in the graph C*-algebra case is that our techniques are much simpler than Katsura's construction of the topological graph algebra and its deep structure results. Regarding the algebraic setting we are unaware of a converse for the Cuntz-Krieger uniqueness theorem for Leavitt path algebras and believe this is a new result.

This paper is organized as follows: In Section~2 we give a review necessary to make the paper self-contained. In Section~3 we present a sufficient condition over branching systems of a directed graph such that the representation induced from a branching system satisfying this condition is faithful. Then we construct a class of branching systems associated to a directed graph satisfying the above condition. This class of examples were firstly built in \cite{GR3}, in the algebraic setting, and hence it is interesting to note that the same class of branching systems provide faithful representations of both Leavitt path algebras and graph C*-algebras. We finish this paper by proving the converse of the Cuntz-Krieger uniqueness theorem for graph algebras.

\section{Preliminaries}

Throughout this paper, all measure spaces are assumed to be $\sigma$-finite.

In this section we recall some background about directed graphs and its corresponding algebras. We also recall the notion of branching systems of a directed graph and the construction of a representation of the graph algebra from a branching system.

Firstly recall that a \emph{directed graph} is a quadruple $E=(E^0, E^1, r,s)$ consisting of two countable sets $E^0, E^1$, and two maps $r,s:E^1 \to E^0$. We think of $E^0$ as a set of vertices, and we think of every element $e \in E^1$ as an arrow pointing from $s(e)$ to $r(e)$. The graph $E$ is called \emph{row-finite} if $\vert s^{-1}(v)\vert< \infty$ for all $v \in E^0$. For $v \in E^0$, we call $v$ a \emph{sink} if $s^{-1}(v)=\emptyset$, and we call $v$ a \emph{source} if $r^{-1}(v)=\emptyset$. In this paper we use the following combinatorial definitions regarding a directed graph:

\begin{definicao}[{\cite{Raeburn:Graphalgebras05, Szyma'nski:IJM02}}]\label{defofConditionLforultragraph}
Let $E$ be a directed graph. For $n \geq 1$, a \emph{path} of length $n$ is a tuple $(e_i)_{i \in 1}^{n} \in \prod_{i=1}^{n}E^1$ such that $r(e_i)=s(e_{i+1})$ for $i=1,\dots,n-1$. The path $(e_i)_{i=1}^{n}$ is called a \emph{cycle} if $s(e_1) =r(e_n)$, and $s(e_1)$ is called the \emph{base point} of the cycle. The cycle is called \emph{simple} if $r(e_i)\neq r(e_j)$ for all $i \neq j$. We say the cycle $(e_i)_{i=1}^{n}$ has no exits if $r^{-1}(r(e_i))=e_i$ for all $i$. We say that the graph $E$ satisfies \emph{Condition~(L)} if any cycle of $E$ has no exits.
\end{definicao}

Recall that the \emph{graph C*-algebra} $C^*(E)$, as defined in \cite{flr}, is the universal C*-algebra generated by a family of partial isometries with orthogonal ranges $\{s_e:e\in E^1\}$ and a family of mutually orthogonal projections $\{p_v:v\in E^0\}$ satisfying
\begin{enumerate}
\item\label{s_e^*s_e=p_{r(e)}}$s_e^*s_e=p_{r(e)}$, for all $e\in E^1$;
\item $s_es_e^*\leq p_{s(e)}$ for all $e\in E^1$; and
\item\label{CK-condition} $p_v=\sum\limits_{s(e)=v}s_es_e^*$ whenever $0<\vert s^{-1}(v)\vert< \infty$.
\end{enumerate}

Leavitt path algebras may be defined in terms of the same relations as above, though in the algebraic context the more common definition is the following one: Given a graph $E$ and a field $K$, the \emph{Leavitt path algebra} of $E$, denoted by $L_K(E)$, is the universal $K$-algebra generated by a set $\{v: v \in E^0\}$, of pairwise orthogonal idempotents, together with a set $\{e, e^* : e \in E^1\}$ of elements satisfying  
\begin{enumerate}
\item $s(e)e = er(e) = e$, $r(e)e^* = e^*s(e) = e^*$ and $e^*f = \delta_{e,f} r(e)$ for all $e,f \in E^1$,
\item $v =\sum\limits_{e\in E^1:s(e)=v}ee^*$ for every vertex $v$ with $0 < \#\{e: s(e) = v\} <\infty.$
\end{enumerate}

Now we recall the notion of branching systems of a directed graph from \cite{MR2903145}.

\begin{definicao}[{\cite[Definition~2.1]{MR2903145}}]\label{branchsystem}
Let $E$ be a directed graph, let $(X,\mu)$ be a measure space, and let $\{R_e,D_v\}_{e\in E^1,v\in E^0}$ be a family of measurable subsets of $X$. Suppose that
\begin{enumerate}
\item $R_e\cap R_f \stackrel{\mu-a.e.}{=}\emptyset$ if $e \neq f \in E^1$;
\item $D_v \cap D_w\stackrel{\mu-a.e.}{=}\emptyset$ if $v \neq w \in E^0$;
\item $R_e\stackrel{\mu-a.e.}{\subseteq}D_{s(e)}$ for all $e\in E^1$;
\item $D_v\stackrel{\mu-a.e.}{=} \bigcup_{e \in s^{-1}(v)}R_e$ if $0 <\vert s^{-1}(v) \vert<\infty$; and
\item for each $e\in E^1$, there exist two measurable maps $f_e:D_{r(e)}\rightarrow R_e$ and $f_e^{-1}:R_e \rightarrow D_{r(e)}$ such that $f_e\circ f_e^{-1}\stackrel{\mu-a.e.}{=}\id_{R_e}, f_e^{-1}\circ f_e\stackrel{\mu-a.e.}{=}\id_{D_{r(e)}}$, the pushforward measure $\mu \circ f_e$ of $f_e^{-1}$ in $D_{r(e)}$ is absolutely continuous with respect to $\mu$ in $D_{r(e)}$, and the pushforward measure $\mu \circ f_e^{-1}$ of $f_e$ in $R_e$ is absolutely continuous with respect to $\mu$ in $R_e$. Denote the Radon-Nikodym derivative $d(\mu \circ f_e)/d\mu$ by $\Phi_{f_e}$, and denote the Radon-Nikodym derivative $d(\mu\circ f_e^{-1} )/d\mu$ by $\Phi_{f_e^{-1}}$.
\end{enumerate}
We call $\{R_e,D_v,f_e\}_{e \in E^1,v \in E^0}$ an \emph{$E$-branching system} on the measure space $(X,\mu)$.
\end{definicao}

\begin{obs}
In the algebraic context, an $E$-algebraic branching system as defined in \cite{GR3} is the same as an $E$-branching system, except we deal with exact equalities instead of equality almost everywhere, there is no mention to measures or to Radon-Nykodym derivatives and the maps between the sets are only required to be bijections.
\end{obs}

\begin{teorema}[{\cite[Theorem~2.2]{MR2903145}}]\label{repinducedbybranchingsystems} 
Let $E$ be a directed graph. Fix an $E$-branching system $\{R_e,D_v,f_e\}_{e\in E^1, v\in E^0}$ on a measure space $(X,\mu)$. Then there exists a unique representation $\pi:C^*(E) \to B(\mathcal{L}^2(X,\mu))$ such that $\pi(s_e)(\phi)=\Phi_{f_e^{-1}}^{1/2}( \phi \circ f_e^{-1})$ and $\pi(p_v)(\phi)=\chi_{D_v}\phi$, for all $e \in E^1, v \in E^0$, and for all $\phi \in \mathcal{L}^2(X,\mu)$.
\end{teorema}

\begin{obs}
In a similar way as above, see \cite{GR3}, given an $E$-algebraic branching system we obtain a representation $\pi$ of $L_K(E)$ in $\text{Hom}_K(M)$, the $K$ algebra of all homomorphism from $M$ to $M$, where $M$ is the $K$ module of all functions in $X$, such that for all $v\in E^0$, $e\in E^1$, and $\phi \in M$, $\pi(v)(\phi)= \chi_{D_v} \phi$, $\pi(e)(\phi) = \chi_{R_e} \cdot \phi \circ f_{e^{-1}}$ and $\pi(e^*)(\phi) = \chi_{D_{r(e)}} \cdot \phi \circ f_{e}$. 
\end{obs}

Finally, for an $E$-branching system $\{R_e,D_v,f_e\}_{e\in E^1,v\in E^0}$ on a measure space $(X,\mu)$, let $\pi:C^*(E) \to B(L^2(X,\mu))$ be the representation induced from the branching system. Fix a finite path $\alpha \in E^n$, for some $n \geq 1$. Define $f_\alpha:=f_{\alpha_1} \circ\dots\circ f_{\alpha_n}$, and define $f_\alpha^{-1}:=f_{\alpha_n}^{-1} \circ\dots\circ f_{\alpha_1}^{-1}$. It is straightforward to see that $\mu \circ f_{\alpha_1} \circ\dots\circ f_{\alpha_n}$ in $D_{r(\alpha_n)}$ is absolutely continuous with respect to $\mu$ in $D_{r(\alpha_n)}$, and $\mu \circ f_{\alpha_n}^{-1} \circ\dots\circ f_{\alpha_1}^{-1}$ in $R_{\alpha_1}$ is absolutely continuous with respect to $\mu$ in $R_{\alpha_1}$. Denote the Radon-Nikodym derivative $d(\mu \circ f_{\alpha_1}\circ\dots\circ f_{\alpha_n})/d\mu$ by $\Phi_{f_\alpha}$, and denote the Radon-Nikodym derivative $d(\mu \circ f_{\alpha_n}^{-1}\circ\dots\circ f_{\alpha_1}^{-1})/d\mu$ by $\Phi_{f_\alpha^{-1}}$. So for any $\phi \in L^2(X,\mu)$, we have that
\begin{equation}
\pi(s_\alpha)(\phi)=\Phi_{f_{\alpha}^{-1}}^{1/2}\phi \circ f_{\alpha}^{-1}, \text{ and } \pi(s_\alpha)^*(\phi)=\Phi_{f_{\alpha}}^{1/2}\phi \circ f_{\alpha}
\end{equation}
and the analogue result also holds in the algebraic context.

\section{Faithful Representations}

For row-finite directed graphs without sinks in \cite[Theorem~4.2]{GR3} it was shown that, under a mild condition over an algebraic branching system, the induced Leavitt path algebra representation is faithful. Next, for any directed graph $E$, we give an analogous condition over an $E$-branching system so that the induced graph C*-algebra representation is faithful. The following theorem is our main result in this paper.

\begin{teorema}\label{a criterion of faithful rep}
Let $E$ be a directed graph, $\{R_e,D_v,f_e\}_{e\in E^1,v\in E^0}$ be an $E$-branching system on a measure space $(X,\mu)$ such that $\mu(D_v) \neq 0$ for all $v \in E^0$ and let $\pi:C^*(E) \to B(L^2(X,\mu))$ be the representation induced from the branching system. Suppose that for each $v \in E^0$ such that $v$ is a base point of a cycle which has no exits, and for finitely many cycles $\{\alpha^i\}_{i=1}^{n}$ with the base point $v$, there exists a measurable subset $F$ of $D_v$ with $\mu(F) \neq 0$, such that $f_{\alpha^i}(F) \cap F\stackrel{\mu-a.e.}{=}\emptyset$ for all $i$. Then $\pi$ is faithful.
\end{teorema}

\demo
For each $v \in E^0$, since $\mu(D_v) \neq 0$, we have that $\pi(p_v) \neq 0$. For $v \in E^0$ such that $v$ is a base point of a cycle which has no exits, there exists a unique simple cycle $\alpha=(e_1,...,e_m)$ with the base point $v$. In order to show that $\pi$ is faithful, by \cite[Theorem~1.2]{Szyma'nski:IJM02}, we only need to show that the spectrum of $\pi(s_\alpha)$ contains the full circle. Since $\alpha$ is a simple cycle without exits, by the universal property of $C(\mathbb{T})$, there exists a unique homomorphism $h:C(\mathbb{T}) \to C^*(\pi(s_{\alpha}))$, such that $h(I)=\pi(p_v)$, and $h(u)=\pi(s_\alpha)$ where $u$ is the universal unitary element in $C(\mathbb{T})$. Since the spectrum of $u$ is the full circle, by \cite[Corollary~II.1.6.7]{BL}, to prove that the spectrum of $\pi(s_\alpha)$ in $C^*(E)$ contains the full circle, it is sufficient to prove that $h$ is an isomorphism. By \cite[Proposition~3.2]{Raeburn:Graphalgebras05} (by considering the action $\mathbb{T}\ni z\mapsto \beta_z\in Aut(C(\mathbb{T}))$ defined by $\beta_z(u)=zu$ and $\beta_z(I)=I$), there exists a faithful conditional expectation on $C(\mathbb{T})$ sending $u^i (u^* )^j$ to $\delta_{i,j}I$. Hence by \cite[Proposition~3.11]{Katsura:CJM03}, to show that $h$ is faithful, we only need to construct a conditional expectation $\Psi:C^*(\pi(s_{\alpha})) \to C^*(\pi(s_{\alpha}))$ such that $\Psi(\pi(s_{\nu})\pi(s_\tau)^*)=\delta_{\vert\nu\vert,\vert\tau\vert}\pi(s_{\nu})\pi(s_\tau)^*$ whenever $s(\nu)=s(\tau)=r(\nu)=r(\tau)=v$. Since 
\[
C^*(\pi(s_{\alpha})):=\overline{\mathrm{span}}\{\pi(p_v),\pi(s_\nu),\pi(s_\tau)^*:s(\nu)=s(\tau)=r(\nu)=r(\tau)=v\},
\]
it is sufficient to show that 

$$\vert z\vert\leq\Big\Vert z\pi(p_v)+\sum_{i=1}^{n}z_i\pi(s_{\nu^i})+\sum_{j=1}^{m}z_j'\pi(s_{\tau^j})^*\Big\Vert$$
whenever $s(\nu^i)=r(\nu^i)=s(\tau^j)=r(\tau^j)=v$ for each $i,j$ and $z,z_i, z_j' \in \C$. We do this below. 

By the assumption of the theorem, there exists a measurable subset $F$ of $D_v$ with $\mu(F) \neq 0$, such that $f_{\nu^i}(F) \cap F\stackrel{\mu-a.e.}{=}\emptyset$ and $f_{\tau^j}(F) \cap F\stackrel{\mu-a.e.}{=}\emptyset$ for each $i,j$. Take an arbitrary function $\phi \in L^2(X,\mu)$ with $\Vert\phi\Vert = 1$ and $\supp(\phi) \stackrel{\mu-a.e.}{\subset}F$. Then $\pi(s_{\nu^i})(\phi)(x)=0$ and $\pi(s_{\tau^j})^*(\phi)(x)=0$ for each $i,j$ and almost every $x\in F$.
Then
\begin{align*}
\Big\Vert z\pi(p_v)(\phi)&+\sum_{i=1}^{n}z_i\pi(s_{\nu^i})(\phi)+\sum_{j=1}^{m}z_j'\pi(s_{\tau^j})^*(\phi)\Big\Vert^2
\\&=\int_X  \Big\vert z\pi(p_v)(\phi)+\sum_{i=1}^{n}z_i\pi(s_{\nu^i})(\phi)+\sum_{j=1}^{m}z_j'\pi(s_{\tau^j})^*(\phi)\Big\vert^2 \, \mathrm{d}\mu
\\&\geq\int_F  \Big\vert z\pi(p_v)(\phi)+\sum_{i=1}^{n}z_i\pi(s_{\nu^i})(\phi)+\sum_{j=1}^{m}z_j'\pi(s_{\tau^j})^*(\phi)\Big\vert^2 \, \mathrm{d}\mu
\\& = \int_F\Big\vert z\pi(p_v)(\phi)\Big\vert^2\,\mathrm{d}\mu 
\\&=\vert z\vert^2.
\end{align*}
So $\vert z \vert\leq\Big\Vert z\pi(p_v)+\sum_{i=1}^{n}z_i\pi(s_{\nu^i})+\sum_{j=1}^{m}z_j'\pi(s_{\tau^j})^*\Big\Vert$ and hence we are done.

\fim


Next we introduce a class of branching systems satisfying the condition of the previous theorem.

Let $E$ be a directed graph, let $X=\mathbb{R}$ and let $\mu$ be the Lebesgue measure on all Borel sets of $\mathbb{R}$. Enumerate $E^1=\{e_i\}_{i \geq 1}$ and the set of sinks $E_{\mathrm{sink}}^0=\{v_i:s^{-1}(v_i)=\emptyset\}_{i \geq 1}$, where each $i$ is a natural number. For each $i \geq 1$, define $R_{e_i}:=[i-1,i)$ and $D_{v_i}:=[-i,1-i)$. For $v \in E^0$, with $s^{-1}(v) \neq \emptyset$, define $D_v:=\cup_{e \in s^{-1}(v)}R_e$. Now, for each $e \in E^1$, define $f_e$ as an arbitrary diffeomorphism $f_e:D_{r(e)} \to R_e$ and denote the derivative of $f_e$ by $\Phi_{f_e}$ and the derivative of $f_e^{-1}$ by $\Phi_{f_e^{-1}}$. By \cite[Theorem~3.1]{MR2903145}, we have that $\{R_e,D_v,f_e\}_{e \in E^1,v \in E^0}$ is an $E$-branching system on $(X,\mu)$. Let $\pi:C^*(E) \to B(L^2(X,\mu))$ be the induced representation.

In the next paragraph we redefine some of the maps $f_e$ defined above to obtain branching systems that induce faithful representations of $C^*(E)$.

Denote by $W$ the set of vertices which are base points of cycles without exits. For each $w \in W$, there exists a unique simple cycle $\alpha=(\alpha_1,...,\alpha_m)$ with the base point $w$. Notice that $D_{r(\alpha_i)}$ and $R_{\alpha_i}$ are all unit intervals, that is, for $1 \leq i \leq m,D_{r(\alpha_i)}=[k_i,k_i+1)$ and $R_{\alpha_i}=[l_i,l_i+1)$, for some $k_i,l_i \geq 0$. For $1 \leq i \leq m$, take $\theta_i \in [0,1)$ and define $f_{\alpha_i}:D_{r(\alpha_i)}\rightarrow R_{\alpha_i}$ by $f_{\alpha_i}(x)=(x+\theta_i )mod(1)+l_i$ (instead of any diffeomorphpism). So we now have a new  $E$-branching system and below we characterize when this branching system induces a faithful representation of $C^*(E)$.

For each $w\in W$, consider the unique simple cycle $\alpha=(\alpha_1,...,\alpha_m)$ whose base point is $w$ and let $f_\alpha:D_{r(\alpha_m)}\rightarrow R_{\alpha_1}$ be the composition $f_\alpha=f_{\alpha_1}\circ...\circ f_{\alpha_m}$. Since $D_{r(\alpha_m)}=R_{\alpha_1}=[l_1,l_1+1)$ we get that $f_\alpha:[l_1,l_1+1)\rightarrow [l_1,l_1+1)$.

It is not hard to see (by direct calculations) that $$f_\alpha(x)=[x+(\theta_1+\theta_2+...+\theta_m)mod(1)]mod(1)+l_1,$$ for each $x\in [l_1,l_1+1)$. Let $\theta_w=(\theta_1+\theta_2+...+\theta_m)mod(1)$ and notice that $f_\alpha(x)=(x+\theta_w)mod(1)+l_1$ for each $x\in [l_1,l_1+1)$.

\centerline{
\setlength{\unitlength}{1.5cm}
\begin{picture}(3,3.5)
\put(-0.5,0){\vector(1,0){5}}
\put(0.5,-0.5){\vector(0,1){3,5}}
\put(1.5,-0.05){\line(0,1){0.1}}
\put(1.4,-0.3){$l_1$}
\put(3.5,-0.05){\line(0,1){0.1}}
\put(3.3,-0.3){$l_1+1$}
\put(0.45,1){\line(1,0){0.1}}
\put(0.2,0.9){$l_1$}
\put(0.45,2.5){\line(1,0){0.1}}
\put(-0.2,2.4){$l_1+1$}
\put(0.45,2){\line(1,0){0.1}}
\put(-0.35,1.85){$l_1+\theta_w$}
\put(1.5,2){\line(1,1){0.46}}
\put(2,2.5){\circle{0.1}}
\put(1.5,2){\circle*{0.1}}
\put(2,-0.05){\line(0,1){0.1}}
\put(1.8,-0.3){$l_1+1-\theta_w$}
\put(2,1){\line(1,1){1.47}}
\put(2,1){\circle*{0.1}}
\put(3.5,2.5){\circle{0.1}}
\put(1.5,-1){Graph of $f_{\alpha}$}
\end{picture}}
\vspace{2cm}


\begin{proposicao} Let $\{R_e,D_v,f_e\}_{e\in E^1, v\in E^0}$ be the branching system introduced above and let $\pi:C^*(E)\rightarrow B(L^{2}(X,\mu))$ be the induced representation. Then $\pi$ is faithful if and only if $\theta_w$ is irrational for each $w\in W$.
\end{proposicao}

\demo
First suppose that each $\theta_w$ is irrational. By Theorem~\ref{a criterion of faithful rep} it is enough to show that, for finitely many cycles $\{\beta^i\}_{i=1}^{n}$ with the base point $w$, there exists a measurable subset $F$ of $D_w$, with $\mu(F) \neq 0$, such that $f_{\beta^i}(F) \cap F\stackrel{\mu-a.e.}{=}\emptyset$ for all $i$. 

Notice that each $\beta^i$ has the form $\beta_i=(\alpha,\dots,\alpha)$ ($q_i$ times), where $\alpha$ is the unique simple cycle based on $w$. By direct calculations it follows that $f_{\beta_i}(x)=(x+(q_i\theta_w))mod(1)+l_1$, for each $x\in D_w=[l_1,l_1+1)$ and hence (looking at the graph of $f_{\beta_i}$) we have that $f_{\beta_i}([l_1,y))=[f_{\beta_i}(l_1), f(y))$, for each $y\in [l_1,l_1+1-(q_i\theta_w)mod(1))$. Since $f_{\beta_i}(l_1)=l_1+(q_i\theta_w)mod(1)$ and $\theta_w$ is irrational then $f_{\beta_i}(l_1)$ is irrational and so $f_{\beta_i}(l_1)> l_1$ for each 
$\beta_i$. 

Now, let $c\in \R$ be such that $c>l_1$, $c<l_1+1-(q_i\theta_w)mod(1)$ and $c<f_{\beta_i}(l_1)$, for each $i\in \{1,...,n\}$, and define $F=[l_1,c)$. Then $\mu(F)\neq 0$ and $f_{\beta_i}(F)\cap F=\emptyset$ for each $\beta_i$ and hence, by Theorem~\ref{a criterion of faithful rep}, we have that $\pi$ is faithful.

Suppose now that some $\theta_w\in [0,1)$ is rational, say $\theta_w=\frac{p}{q}$ with $p,q$ positive integers. Let $\alpha$ be the (unique) simple cycle based on $w$ and let $\beta=(\alpha,...,\alpha)$ ($p$ times). Note that for each $x\in D_w=[l_1,l_1+1)$ we have that $$f_\beta(x)=[x+(p\theta_w)mod(1)]mod(1)+l_1=(x)mod(1)+l_1=x,$$ and therefore $\pi(S_\beta)=\pi(p_w)$ and $\pi$ is not faithful. 

\fim

\begin{obs}
The above result allows us to construct faithful representations of graph C*-algebras even when the condition of the Cuntz-Krieger uniqueness theorem fails. We exemplify below.
\end{obs}

\begin{exemplo}\label{example of a graph with a single cycle}
Let $E$ be a row finite directed graph consisting of a single cycle of length $1$, that is, $E^0=\{v\}$,  $E^1=\{e\}$, $r(e)=s(e)=v$. Let $X=\mathbb{R}$ and let $\mu$ be the Lebesgue measure on all Borel sets of $\mathbb{R}$. Fix an irrational number $\theta \in [0,1)$. Define $D_v=R_e:=[0,1)$, and define $f_e:D_{v} \to R_e$ by $f_e(x):=(x+\theta)mod(1)$. Then $\{R_e,D_v,f_e\}$ is an $E$-branching system. By the above discussions, the representation induced by this branching system is faithful.
\end{exemplo}

We mention that Katsura proved a version of the converse of the Cuntz-Krieger uniqueness theorem for topological graph algebras (see \cite[Theorem~6.14]{Katsura:ETDS06}), whose proof is very complicated. The following theorem is an application of branching systems which gives a simple proof of the converse of the Cuntz-Krieger uniqueness theorem for graph algebras.

\begin{teorema}\label{a criterion of nonfaithful rep}
Let $E$ be a directed graph not satisfying Condition~(L). Then there exist an $E$-branching system $\{R_e,D_v,f_e\}$ on a measure space $(X,\mu)$ and a representation $\pi:C^*(E) \to B(L^2(X,\mu))$ from Theorem~\ref{repinducedbybranchingsystems} such that $\pi(p_v)\neq 0$ for all $v \in E^0$ and $\pi$ is not faithful.
\end{teorema}

\demo
Since $E$ does not satisfy Condition~(L), there is a cycle $\alpha=(\alpha_1,\dots,\alpha_n)$ such that $\alpha_i \neq \alpha_j$ if $i \neq j$, and $s^{-1}(s(\alpha_i))=\{\alpha_i\}$ for all $i$. We enumerate the edge set as $E^1=\{\alpha_1,\dots,\alpha_n,e_{n+1},\dots\}$, and enumerate the vertex set as $E^0=\{s(\alpha_1),\dots,s(\alpha_n),v_{n+1},\dots\}$. By the construction in \cite[Theorem~3.1]{MR2903145}, there is an $E$-branching system on $(\mathbb{R},\mu)$ denoted by $\{R_e,D_v,f_e\}$, where $\mu$ is the Lebesgue measure on all Borel sets of $\mathbb{R}$, such that for each $i, \ D_{s(\alpha_i)}=R_{\alpha_i}=[i-1,i]$ and $f_{\alpha_i}$ is the increasing bijective linear map. Notice that $f_\alpha=\id$ and so $\Phi_{f_\alpha}\equiv 1$ on $R_{\alpha_1}=[0,1]$. So $\pi(s_\alpha^*)=\pi(p_{s(\alpha_1)})$. By the construction in \cite[Theorem~3.1]{MR2903145}, we deduce that $\pi(p_v)\neq 0$ for all $v \in E^0$.

Suppose that $\pi$ is not faithful, for a contradiction. By the universal property there exists a gauge action $\gamma$ on $\pi(C^*(E))$. So for each $z \in \mathbb{T}$ we have that
\[
\pi(p_{s(\alpha_1)})=\gamma_z(\pi(p_{s_{\alpha_1}}))=\gamma_z(\pi(s_\alpha^*))=\overline{z}^n\pi(s_\alpha^*)=\overline{z}^n\pi(p_{s(\alpha_1)}),
\]
which is impossible. Therefore $\pi$ is not faithful.
\fim

Using the theory of branching systems we can also prove the converse of the Cuntz-Krieger uniqueness theorem (see \cite{Tomforde}) for Leavitt path algebras, a result that we could not find in the literature.

\begin{teorema}\label{Leavitt converse CKU}
Let $E$ be a directed graph not satisfying Condition~(L). Then there exists an $E$-algebraic branching system $\{R_e,D_v,f_e\}$ such that the representation $\pi:L_K(E) \to Hom (M)$ given above do not vanishes at the vertices, that is $\pi(v) \neq 0$ for all $v\in E^0$, but $\pi$ is not faithful.
\end{teorema}

\demo
Since $E$ does not satisfy Condition~(L), there exists a cycle $\alpha=(\alpha_1,\dots,\alpha_n)$ such that $\alpha_i \neq \alpha_j$ if $i \neq j$, and $s^{-1}(s(\alpha_i))=\{\alpha_i\}$ for all $i$. 

We enumerate the edge set $E^1=\{\alpha_1,\dots,\alpha_n,e_{n+1},\dots\}$, and enumerate the vertex set as $E^0=\{s(\alpha_1),\dots,s(\alpha_n),v_{n+1},\dots\}$. Following the construction given in theorem 3.1 of \cite{GR3} (which is analogous to the construction presented above for graph C*-algebras), we obtain an $E$-algebraic branching system on $\mathbb{R}$, such that for each $i, \ D_{s(\alpha_i)}=R_{\alpha_i}=[i-1,i)$ and $f_{\alpha_i}$ is the increasing bijective linear map. So $\pi(s_\alpha^*)=\pi(p_{s(\alpha_1)})$. 

To complete the proof we need to show that $s_\alpha^*\neq p_{s(\alpha_1)}$ in $L_K(E)$. But this can be done using once more the theory of algebraic branching systems. Just notice that, if in the construction above instead of picking $f_{\alpha_1}: D_{r(\alpha_1)} \rightarrow R_{\alpha_1}$ as the increasing bijective linear map we pick $f_{\alpha_1}$ as a non-linear bijective increasing map, and we keep the same choice for the remaining $f_{\alpha_i}$, then $f_\alpha \neq id$ and it is straightforward to check that $\pi(s_\alpha^*-p_{s(\alpha_1)}) \neq 0 $ and hence $s_\alpha^*\neq p_{s(\alpha_1)}$  as desired.
\fim

\section*{Acknowledgments}

The second author would like to thank Australian Mathematical Society for offering him a Lift-off Fellowship, he would like to thank Professor David Pask and Professor Aidan Sims for recommending him for the consideration of a Lift-off Fellowship, and he would like to thank Dr Ngamta Thamwattana for providing useful information about the application of the fellowship. The second author in particular wants to thank the extreme hospitality of Professor Daniel Gon\c{c}alves, Professor Danilo Royer and Departamento de Matem\'{a}tica, Universidade Federal de Santa Catarina during his stay in Brazil. Finally the second author appreciates the support from Research Center for Operator Algebras of East China Normal University.

\vspace{1.5pc}
\begin{center}
Daniel Gon\c{c}alves (daemig@gmail.com) and Danilo Royer (danilo.royer@ufsc.br)\\

Departamento de Matem\'{a}tica - Universidade Federal de Santa Catarina, Florian\'{o}polis, 88040-900, Brazil

\vspace{1.5pc}

Hui Li (hli@math.ecnu.edu.cn)\\

Research Center for Operator Algebras, Department of Mathematics, East China Normal University (Minhang Campus), 500 Dongchuan Road, Minhang District, Shanghai 200241, China
\end{center}

\end{document}